\documentclass[12pt]{amsart}
\usepackage{latexsym,amsfonts,amsthm,amsmath,amscd,amssymb}

\newtheorem{teo}{\bf Theorem}
 
 \newtheorem{lem}[teo]{\bf Lemma}
 \newtheorem{prop}[teo]{\bf Proposition}
\newtheorem{opr}[teo]{\bf Definition}
\newcommand {\vse} {$\blacksquare$}
\newcommand \dist {\mathop{\rm dist}\nolimits}
\newcommand{\0}{\varnothing}

\begin{document}
\author{Vladislav Aseev}
\address{Sobolev Institute of Mathematics\\ Novosibirsk\\
Russia}\email{ase@math.nsc.ru}

\author{Tatyana
Kergylova}
\address{\noindent Department of Mathematics\\
Gorno-Altaisk state University \\ Russia}\email{kergyl@gmail.com}

\title[Moebius transformations preserving fixed anharmonic ratio]{Moebius transformations \\ preserving fixed \\ anharmonic ratio}

\thanks{Authors were supported
by the grant of the RFBR N 07-01-90800-mob-st}

\begin{abstract}
O. Kobayashi \cite{kob} proved that $C^1$-mappings preserving
anharmonic ratio are Moebius transformations. We strengthen his
result and prove, that the requirement of differentiability and
even of injectivity can be omitted.
\end{abstract}

\date{\today}

\maketitle

\section{Introduction}

 A concept of Apollonian tetrad in the complex plane
$\textbf{C}$ was introduced in the paper (\cite{haruki},\ Def. 1,\
p.15): any ordered quadruple of distinct points $\{z_1,\ z_2,\
z_3,\ z_4\}\subset\textbf{C}$ is called an {\em Apollonian}
tetrad, if
$$
|z_2-z_3|\cdot|z_1-z_4|=|z_3-z_1|\cdot|z_2-z_4|=|z_1-z_2|\cdot|z_3-z_4|
\eqno (0.1) $$ It was proved in (\cite{haruki},\ Main Theorem,\ p.
19), that any univalent analytic function in the domain
$D\subset\textbf{C}$ is linear-fractional iff the image of any
Apollonian tetrad in $D$ is also Apollonian tetrad.

Since that time more than ten articles appeared, in which
different generalizations of the mentioned property of tetrad to
finite set of points were introduced:  Haruki~H. and Rassias~ T.M.
in \cite{hh} considered Apollonian triangles and hexagons, Bulut
S. and Yilmaz \"Ozg\"ur N. in  \cite{bul} considered Apollonian
set consisting of $2 n$ pairwise different points and proved, that
the analytic univalent function is linear-fractional iff the image
of any Apollonian set in $D$ is also an Apollonian set.

The constructions in these works have considerable interest from
the point of view of the theory of quasimoebius transformations in
$\overline{\bf{R}}^n$ (see \cite{va}) and in Ptolemaic spaces (see
\cite{as}).

O. Kobayashi \cite{kob} noticed, that relation~(0.1) is equivalent
to the equality $[z_1: z_2:z_3: z_4]=(1\pm i \sqrt{3})/2$, where
$[z_1: z_2:z_3: z_4]$ is the  anharmonic ratio of the quadruple
$\{z_1,\ z_2,\ z_3,\ z_4\}$, and he obtained the following result:

\begin{teo}\label{t1} {\bf( \cite{kob},\ Theorem 2.1,\ p. 118)}

Let $\lambda\in \textbf{C}\setminus \textbf{R}$  and
$U\subset\textbf{C}$ be a domain. Suppose
$f:U\rightarrow\textbf{C}$ is an injective differentiable mapping.

If for any quadruple of pairwise distinct points $\{z_1,\ z_2,\
z_3,\ z_4\}\in~
 U$ with anharmonic ratio $[z_1: z_2:z_3: z_4]=\lambda$ the
 equality
  $[f(z_1): f(z_2):f(z_3): f(z_4)]=\lambda$ holds,
 then $f$ is a Moebius transformation.
 \end{teo}

In this paper we show, that the theorem is valid  for any $\lambda
\in \overline{\textbf{C}}\setminus \{0,\ 1,\ \infty\}$  and for
any continuous non-constant mapping
$f:U\rightarrow\overline{\textbf{C}}$ of the domain
$U\subset\overline{\textbf{C}}$ without the requirement of
injectivity and differentiability of $f$.

\section{The tables of symbols,  terminology and the main result }
Let $\textbf{T}$ be set of all ordered quadruples (tetrads)
$T=\{z_1,\ z_2,\ z_3,\ z_4\}$ in  the extended complex plane
$\overline{\textbf{C}}$, that do not contain any three coincident
elements.

A tetrad with four pairwise different elements is called {\it
nonsingular}.

For any tetrad we define its anharmonic ratio (see \cite{berdon},\
\S 44)
 $A(T)=[z_1:
z_2:z_3:~ z_4~].$

If all points in a tetrad $T$ are finite and pairwise different,
then $A(T)$ is
$$A(T)=\dfrac{(z_1-z_3)(z_2-z_4)}{(z_3-z_2)(z_4-z_1)}.\eqno(1.1.1)$$

Under these conditions $A(T)$ is different from $0,\ 1$ and
$\infty$. For tetrads $s_{12}(T),\ s_{13}(T),\ s_{14}(T)$,
obtained from $T$ by permutation of first and  second, of first
and third, of first and fourth elements respectively the following
equalities hold:
$$A(s_{12}(T))=\dfrac{1}{A(T)};\ A(s_{13}(T))=\dfrac{A(T)}{A(T)-1};\
A(s_{14}(T))=1-A(T)\eqno(1.1.2)$$

If in the tetrad $T=\{z_1,\ z_2,\ z_3,\ z_4\}$ neither three
elements coincide, there exists a finite or infinite limit
$$A(T)=\lim\limits_{w_1\rightarrow z_1,\ w_2\rightarrow z_2,\
w_3\rightarrow z_3,\ w_4\rightarrow
z_4}\dfrac{(w_1-w_3)(w_2-w_4)}{(w_3-w_2)(w_4-w_1)},\eqno(1.1.3)$$
where the limit is taken on a set of all nonsingular tetrads
$\{w_1,\ w_2,\ w_3,\ w_4\}$. This limit defines the anharmonic
ratio of the tetrad $T$ in general case. Particularly, for
nonsingular tetrad $T=\{z_1,\ z_2,\ z_3,\ \infty\}$ we have
$A(T)=-(z_1-z_3)/(z_3-z_2)$. We notice, that the tetrad $T\in
\textbf{T}$  is nonsingular iff $A(T)$ different from $0,\ 1$ and
$\infty$.

A Moebius transformation
$\mu:\overline{\textbf{C}}\rightarrow\overline{\textbf{C}}$ is
defined as a superposition of finite number of reflections (or
inversions) with respect to generalized circles in
$\overline{\textbf{C}}$ (see \cite{berdon}, Def 3.1.1, P. 25) and
it realized  either by the linear-fractional function or by its
conjugate.

If $\mu$ is a linear-fractional mapping, then for any tetrad $T\in
\textbf{T}$ the equality $A(\mu(T))=A(T)$ (the invariance of
anharmonic ratio by the linear-fractional mappings) holds. The
opposite is also true:
 a bijective mapping
$\mu:\textbf{C}\rightarrow \textbf{C}$, which preserves an
anharmonic ratio of all nonsingular tetrads, is realized by
linear-fractional function (see \cite{berdon},\ \S4.4).

For given complex number $\alpha\notin \{0,\ 1,\ \infty\}$, we
denote by $\textbf{T}(\alpha)$ the set of all tetrads $T\in
\textbf{T}$ with $A(T)=\alpha$ (all that tetrads are nonsingular)
and by $\textbf{T}(\alpha, \ f)$ the set of all tetrads of
$\textbf{T}(\alpha)$, for which $f(T)\in \textbf{T}$.

\begin{opr} We  say that a mapping $f:U\rightarrow
\overline{\textbf{C}}$, of a domain $U\subset
\overline{\textbf{C}}$  {\em satisfies the condition $\varphi(K,\
\alpha)$}, if for any tetrad $T\in \textbf{T}(\alpha, \ f)$  the
equality $A(f(T))=\alpha=A(T)$ holds.
\end{opr}

 The main result of the paper is the following

\begin{teo}\label{main} If $\alpha \notin \{0,\ 1,\ \infty\}$, then any
continuous mapping $f:U\rightarrow \overline{\textbf{C}}$ of a
domain $U\subset \overline{\textbf{C}}$, which satisfies the
condition $\varphi(K,\ \alpha)$, is either a constant or the
Moebius transformation; and if $\alpha$ is not a real number, then
a function $f(z)$ is either a constant or a linear-fractional
function.\end{teo}

\section{The injectivity lemma }

 \begin{prop}\label{21}
If a continuous mapping $f:U\rightarrow \overline{\textbf{C}}$
satisfies the condition $\varphi(K,\ \alpha)$, $\alpha \notin
\{0,\ 1,\ \infty\}$, then it satisfies the condition $\varphi(K,\
\beta)$, where $\beta$ is taken from the set
$$P_\alpha=\left\{\alpha,\ \dfrac{1}{\alpha},\ 1-\alpha,\
\dfrac{1}{1-\alpha},\ \dfrac{\alpha}{\alpha-1},\
\dfrac{\alpha-1}{\alpha}\right\}\eqno(2.1.1)$$
\end{prop}
{\bf Proof.} The proof follows immediately from relations (1.1.2)
and from the fact, that the equality $f(s(T))=s(f(T))$ is
preserved  by the permutations $s$ of the elements of the tetrad,
and therefore the relations $f(T)\in \textbf{T}$ and $f(s(T))\in
\textbf{T}$ are equivalent.\vse

\begin{lem}\label{inj}
 Suppose $\alpha \notin \{0,\ 1,\ \infty\}$ and let $f:U\rightarrow
 \overline{\textbf{C}}$ be a
continuous mapping  of a domain $U\subset \overline{\textbf{C}}$
satisfying the condition $\varphi(K,\ \alpha)$. Then  $f(z)$ is
either injective in $U$ or $f$ is a constant map.
\end{lem}

{\bf Proof.} We assume that $f$ is not injective. Then we can find
different points $\zeta_0,\ \zeta_\infty\in U$, for which
$f(\zeta_0)=f(\zeta_\infty)=a\in \overline{\textbf{C}}$. Take such
linear-fractional mappings $\mu$ and $\nu$  that
$\mu(\zeta_\infty)=\nu(a)=\infty$. Evidently  $\infty\in
\mu(U)=U^\prime$. A continuous mapping $g=\nu\circ f\circ
\mu^{-1}:U^\prime\rightarrow \overline{\textbf{C}}$ satisfies the
condition $\varphi(K,\ \gamma)$ for any constant $\gamma\in
P_\alpha$; whereas $g(\infty)=g(\mu(\zeta_\infty))=\infty$.

Consider the set
$M=g^{-1}(\{\infty\})\cap\textbf{C}=\{z\in\textbf{C}:g(z)=\infty\}.$

Since $g$ is continuous, the set $M$ is closed in
$U'\backslash\{\infty\}$. We show, that $M$ is an open set.

Let $z_0$  be an any point in $M$.

Build an open disc $B=B(z_0,\ r)\subset U'$ with radius
$$r=\dfrac{\max\{1,\ \dist(z_0,\ \partial U')\}}{2+|\alpha|}\eqno(2.2.1)$$
If $g(z)\equiv \infty$ in $B$, then $B\subset M$, and $z_0$ is an
interior point of the set $M$. Consider the point $z_1\in B$, in
which $g(z_1)\neq\infty$. The point $z_2=z_0+\alpha(z_1-z_0)$ lies
in $U'$, whereas
$$|z_2-z_0|=|\alpha|\cdot|z_1-z_0|\leq \dist(z_0,\ \partial U')\dfrac{|\alpha|}{2+|\alpha|}<\dist(z_0,\ \partial U').$$
For $T=\{z_2,\ z_1,\ z_0,\ \infty\}\subset U'$ we have
$A(T)=-(z_2-z_0)/(z_0-z_1)=\alpha$. The assumption, that
$f(T)=\{g(z_2),\ g(z_1),\ \infty,\ \infty\}\in\textbf{T}$,
contradicts the condition $\varphi(K,\ \alpha)$, because in this
case $A(f(T))\in \{0,\ 1,\ \infty\}$, that is $f(A)\neq \alpha$.

Therefore the tetrad $f(T)$ has three equal elements. But
$g(z_1)\neq\infty$ implies that $g(z_2)=\infty$, that is $z_2\in
M$.

Since $g(z_1)\neq\infty$ and $f$ is continuous in the point $z_1$,
take such $\varepsilon>0$  that the disc $B(z_1,\
\varepsilon)\subset B$ and a function $g(z)\neq \infty$  in
$B(z_1,\ \varepsilon)$. From the inequality $\varepsilon<r$ and
from (2.2.1) it follows that
$$B\left(z_0,\ \dfrac{|\alpha|}{1+|\alpha|}\right)\subset
B\eqno(2.2.2)$$
Putting $\delta=\varepsilon
\dfrac{|\alpha|}{1+|\alpha|}$, consider any point $z\in~ B~(z_0,\
\delta)$. For $w=z+\alpha^{-1} (z_2-z)$ the inequality holds:
$$|w-z_1|=|z-z_1+z_1-z_0+\alpha^{-1}(z_0-z)|<
|z-z_0|\left|1-\dfrac{1}{\alpha}\right|<\varepsilon,$$ which
means, that $w\in B(z_1,\ \varepsilon)$ and therefore
$g(w)\neq\infty$. For $T'=\{z_2,\ w,\ z,\ \infty\}\subset U'$ we
have $A(T')=-(z_2-z)/(z-w)=\alpha$. The assumption that
$g(T')=\{g(z_2)=\infty,\ g(w),\ g(z),\ \infty\}\in\textbf{T}$
contradicts the condition $\varphi(K,\ \alpha)$, because
$A(g(T'))\in \{0,\ 1,\ \infty\}$ and therefore $A(g(T'))\neq
\alpha$. Then in the tetrad $g(T')$ there are three equal
elements. But $g(w)\neq\infty$ and therefore $g(z)=\infty$. Thus
we see, that $g(z)\equiv\infty$ in the disc $B(z_0,\ \delta)$. So
$B(z_0,\ \delta)\subset M$ and $z_0$ is an interior point of the
set $M$. Since any point of the set $M$ it is an interior point,
the set $M$ is open.

Since the set $U'\backslash\{\infty\}$ is connected, the
open-and-closed set $M$ is an empty set or
$U'\backslash\{\infty\}$. But  $\mu(\zeta_0)\in M$, the $M\neq\0$.
Therefore $M=U'\backslash\{\infty\}$, and it means that
$g(z)\equiv\infty$ in $U'$. Thus $f(\zeta)\equiv a$ in $U$.\vse

\section{PROOF OF THE MAIN THEOREM}

We prove the main theorem in several steps arranging them as
independent propositions.

\begin{prop}\label{31}
Suppose the domain $U\subset \overline{\textbf{C}}$ contains
$\infty$, and a continuous injective mapping $f:U\rightarrow
\overline{\textbf{C}},\ f(\infty)=\infty$ satisfies the condition
$\varphi(K,\ \alpha)$, $\alpha \notin \{0,\ 1,\ \infty\}$. Then
for any $z_0\in U$ and $a_0=z_0+w_0,\ b_0=z_0-w_0$, where
$|w_0|<\dist (z_0,\ \partial U)$, the equality holds:
$$f(z_0)=\dfrac{f(a_0)+f(b_0)}{2}\eqno(3.1.1)$$
\end{prop}
{\bf Proof.} Let $\alpha \in \{2,\ 1/2,\ -1\}$. Then by
Proposition \ref{21} the mapping $f$ satisfies the condition
$\varphi(K,\ 2)$. For $T= \{a_0,\ b_0,\ z_0,\ \infty\}$ we obtain
$$2=A(T)=[f(a_0):f(z_0):f(b_0):\infty]=-\dfrac{f(a_0)-f(b_0)}{f(b_0)-f(z_0)},$$
from which the desired equation (3.1.1) immediately follows.

Suppose now $ \alpha \notin \{0,\ 1,\ \infty,\
 2,\ 1/2,\ -1\}$.

If $|\alpha|>1$, then we take $\beta=\alpha$; in otherwise we take
$\beta=1/\alpha$. Then $f$ satisfies the condition $\varphi(K,\
\beta)$, $\beta \notin \{0,\ 1,\ \infty,\ 1/2\}$ (see proposition
\ref{21}). For $\beta'=(1-\beta)/(1-2 \beta)$ we have the equation
$(1-2 \beta)/(1-2 \beta')=-1$. As $\beta\notin \overline{B}(1/2,\
1/2)$, then $|\beta-1/2|>1/2$, that is $|1-2 \beta|>1$. Therefore
$|1-2 \beta'|=q<1$.

Put $R=\dist(z_0,\ \partial U)$ and $w_k={(2 \beta'-1)}^k w_0,\
k=1,\ 2,\ ...$

 As $|w_k|=q^k|w_0|<R$, then all points
$a_k=z_0+w_k$ and $b_k=z_0-w_k$ lies in the disc $B(z_0,\
R)\subset U$. We show that for any $k=0,\
 1,\ 2,\ ...$ the follow equality is valid:
$$f(a_k)+f(b_k)=f(a_0)+f(b_0)\eqno(3.1.2)$$
For $k=0$ (3.1.2) is trivial. We suggest, that it holds for some
$k$ and show, that $f(a_{k+1})+f(b_{k+1})=f(a_0)+f(b_0)$.

For the tetrad $T_1= \{b_k,\ a_{k+1},\ b_{k+1},\ \infty\}$ we have
$$A(T_1)=-\dfrac{b_k-b_{k+1}}{b_{k+1}-a_{k+1}}=-\dfrac{w_{k+1}-w_k}{-2
w_{k+1}}=\dfrac{\beta'-1}{2 \beta'-1}=\beta.$$ The equality
$$A(T_1)=-\dfrac{f(b_k)-f(b_{k+1})}{f(b_{k+1})-f(a_{k+1})}=\beta$$
follows from the condition $\varphi(K,\ \beta)$. For $T_2= \{a_k,\
b_{k+1},\ a_{k+1},\ \infty\}$ we have
$$A(T_2)=-\dfrac{a_k-a_{k+1}}{a_{k+1}-b_{k+1}}=-\dfrac{w_k-w_{k+1}}{2
w_{k+1}}=\beta.$$ Therefore
$$A(T_2)=-\dfrac{f(a_k)-f(a_{k+1})}{f(a_{k+1})-f(b_{k+1})}=\beta.$$
Thus we come to the equality
$$\dfrac{f(b_{k+1})-f(b_k)}{f(b_{k+1})-f(a_{k+1})}=\dfrac{f(a_k)-f(a_{k+1})}{f(b_{k+1})-f(a_{k+1})},$$
from which it follows that $f(b_{k+1})-f(b_k)=f(a_k)-f(a_{k+1})$.

Then by induction we obtain
$$f(a_{k+1})+f(b_{k+1})=f(a_k)+f(b_k)=f(a_0)+f(b_0).$$
Thus we have proved (3.1.2) for all $k=0,\
 1,\ 2,\ ...$

Since
 $w_k\rightarrow 0$, when $k\rightarrow \infty$, then $a_k\rightarrow z_0$ and $b_k\rightarrow
 z_0$.

Next from (3.1.2) we obtain the desired relation (3.1.1). \vse

\begin{prop}\label{32}
Suppose a domain $U\subset \overline{\textbf{C}}$ contains
$\infty$, and a continuous mapping $f:U\rightarrow
\overline{\textbf{C}},\ f(\infty)=\infty$ satisfies the condition
$\varphi(K,\ \alpha)$, $\alpha \notin \{0,\ 1,\ \infty\}$. Then
the mapping $f$ moves any linear segment $L\subset U \backslash
\{\infty\}$ to a linear segment; any ray $P\subset U \backslash
\{\infty\}$ to some ray; any line $Q\subset U \backslash
\{\infty\}$ to some line.
\end{prop}

{\bf Proof.} We show that for any point $c\in U \backslash
\{\infty\}$ the mapping $f$ is linear on any linear segment
$L\subset B(c,\ r)$, where $r=\dist(c,\
\partial U)/3$. Let $a,\ b$ be the endpoints
of the segment $L$. For any point $z\in L$,\ \  $\max\{|z-a|,\
|z-b|\}\leq 2 r<\dist(z,\
\partial U).$  That is, we may apply the Proposition \ref{31} to any point $z\in L$ and  any $w$ such that
$|w|<\max\{|z-a|,\ |z-b|\}$. Therefore for any pair of points
$t_1,\ t_2\in L$ we have $f((t_1+t_2)/2)=(f(t_1)+f(t_2))/2$. We
conclude that the function $f$ is linear  on a dense subset of the
segment $L$ and so, by continuity of $f$, it is linear on $L$.

Thus the mapping $f$ is locally linear on any connected subset
$S\subset U \backslash \{\infty\}$, which lies on a straight line.
Therefore by connectedness of $S$, $f$ is linear on all of $S$.
Particularly, an image of any segment is some segment, the image
of any ray is some ray, the image of any line is some line. \vse

Next we use a criterion of Moebiusness for mappings of $
n$-dimensional domains, proposed by Zelinsky Y.B. for the mapping
of (we take this theorem for case $n=2$).

\begin{teo}\label{ze} {\bf(\cite{zel},\ Th. 8,\ P. 35)}
Suppose a continuous mapping $f:D\rightarrow
\overline{\textbf{C}}$ of a domain $D\subset
\overline{\textbf{C}}$ moves any set $P\subset D$, which lies on a
generalized circle to a set, which lies on a generalized circle.
If $f(D)$ does not lie on a generalized circle, then $f$ is a
Moebius transformation.\end{teo}

\begin{lem}\label{mob} Any continuous injective mapping $f:U\rightarrow
\overline{\textbf{C}}$ of a domain $U\subset
\overline{\textbf{C}}$, which satisfies the condition $\varphi(K,\
\alpha)$,
 $\alpha \notin \{0,\ 1,\
\infty\}$, is a Moebius transformation.
\end{lem}

{\bf Proof.} Take any open disc $D\subset U$ such that $
\overline{D}\subset U$ and a set $P\subset D$ lies on a
generalized circle $S\subset \overline{\textbf{C}}$. Then $S\cap
D$ is a connected subset of a generalized circle $S$.

The situation 1.

Suppose $S\cap \partial D\neq \0$ and $a\in S\cap
\partial D$. We build linear-fractional mappings
$\mu:\overline{\textbf{C}}\rightarrow \overline{\textbf{C}}$ and
$\eta:\overline{\textbf{C}}\rightarrow \overline{\textbf{C}}$ such
that $\mu(a)=\infty,\ \eta(f(a))=\infty$ and consider the mapping
$g=\eta\circ f \circ \mu^{-1}:\mu(U)\rightarrow
\overline{\textbf{C}}$. This mapping satisfies the condition
$\varphi(K,\ \alpha)$ and $\infty\in \mu(U),\ g(\infty)=\infty$.
The set $L'=\mu((S\cap \overline{D})\backslash
\{a\})\subset\mu(U)\backslash \{\infty\}$ is either a ray or a
line. By Proposition \ref{32} the set $g(L')$ is also either some
ray or a line. Therefore $f=\eta^{-1}\circ g \circ \mu$ maps a set
$S\cap D$ (and any subset $P$ also) to  a subset of a circle.

The situation 2.

Take $S\subset D$ and $a\in S$. We build linear-fractional
mappings $\mu$ and $\eta$ by analogy with the situation 1 and
consider the mapping $g=\eta\circ f \circ \mu^{-1}$. By
Proposition \ref{32} the mapping $g$ moves the line
$L'=\mu((S\backslash \{a\}\subset\mu(U)\backslash \{\infty\}$ to
some line. Therefore $g(L')$ is a generalized circle and $f$ moves
the generalized circle $S=\mu^{-1}(L^\prime)$ to a generalized
circle, and any subset $P\subset S$ to a subset of a generalized
circle. Thus $f$ satisfies the conditions of the theorem \ref{ze}
in the disc $D$ and is injective. Therefore $f$ is a Moebius
transformation on the disc $D$. Since our choice of the disc
$D\subset  U$ is arbitrary, the mapping $f$ is locally Moebius on
a domain $U$.  Moebius transformations are explicitly defined by
their values on some quadruple, that does not lie on a generalized
circle. From local Moebiusness of $f$ it follows, that $f$ is
Moebius  on a domain $U$. \vse

\begin{lem}\label{fr} Any continuous injective mapping $f:U\rightarrow
\overline{\textbf{C}}$ of a domain $U\subset
\overline{\textbf{C}}$, which satisfies the condition $\varphi(K,\
\alpha)$, $\alpha \in \textbf{C}\backslash \textbf{R}$ is a
linear-fractional function.
\end{lem}
{\bf Proof.} By Lemma \ref{mob} the mapping $f$ is a Moebius
transformation, so $f(z)$ is either a linear-fractional function
or its conjugate. We show that if $\alpha=a+i b$ and $b\neq 0$,
then the second its impossible. Let $f(z)=\overline{\mu(z)}$,
where $\mu(z)$ is a linear-fractional function. Take any tetrad
$T=\{z_1,\ z_2,\ z_3,\ z_4\}\subset U$ with anharmonic ratio
$A(T)=\alpha=a+i b$. Then $A(\mu(T))=A(T)=a+i b$ and
$A(f(T))=A(\overline{\mu}(T))=\overline{A(\mu(T))}=a-i
b\neq\alpha$. This contradicts the condition $\varphi(K,\
\alpha)$. Therefore $f(z)$ is a linear-fractional function. \vse

The proof of the main Theorem follows immediately from Lemma
\ref{inj}, Lemma \ref{mob} and Lemma \ref{fr}.

\end{document}